\documentclass[10pt]{article}
\usepackage{amsmath, amsthm, amssymb}
\usepackage{authblk}
\usepackage{color}
\usepackage{comment}
\usepackage{enumerate}
\usepackage{float}
\usepackage{graphicx}
\usepackage{hyperref}
\usepackage{cite}

\allowdisplaybreaks[4]
\numberwithin{equation}{section}
\date{}

\theoremstyle{plain}
\newtheorem{theorem}{Theorem}[section]
\newtheorem{proposition}[theorem]{Proposition}
\newtheorem{lemma}[theorem]{Lemma}

\theoremstyle{remark}
\newtheorem{remark}{Remark}
\theoremstyle{definition}

\newcommand{\intd}{\ \mathrm{d}}

\title{A localized criterion for the regularity of solutions to Navier-Stokes equations}

\author{Congming Li$^{a,1}$, Chenkai Liu$^{a,1}$, Ran Zhuo$^{b,c,1}$}

\begin{document}

\renewcommand{\thefootnote}{\fnsymbol{footnote}}

\footnotetext{a. School of Mathematical Sciences and CMA-Shanghai, Shanghai Jiao Tong University, China;}
\footnotetext{b. Mathematics and Science College, Shanghai Normal University, China;}
\footnotetext{c. Department of Mathematical and Statistics, Huanghuai University, China.}

\footnotetext{The authors contributed equally.}
\footnotetext{\textsuperscript{1}To whom correspondence should be addressed. E-mail: congming.li@sjtu.edu.cn; lck0427@sjtu.edu.cn; zhuoran1986@126.com.}


\renewcommand{\thefootnote}{\arabic{footnote}}

\maketitle
\begin{abstract}
The Serrin-Prodi-Ladyzhenskaya type $L^{s,r}$ criteria for the regularity of solutions to the incompressible Navier-Stokes equations are fundamental in the study of the millennium problem posted by the Clay Mathematical Institute about the incompressible  N-S equations.
In this article, we establish some localized $L^{s,r}$ criteria for the regularity of solutions to the equations. In fact, we obtain some a priori estimates of solutions to the equations depend only on some local $L^{s,r}$ type norms. These local $L^{s,r}$ type norms, are small for reasonable initial value and shall remain to be small for global regular solutions. Thus, deriving the smallness or even the boundedness of the local $L^{s,r}$ type norms is necessary and sufficient to affirmatively answer the millennium problem.



\end{abstract}

\maketitle


\section{Introduction}\label{sec1}

Consider the Cauchy problem for the 3-dimensional incompressible Navier-Stokes equations:
\begin{equation}
\label{eq:ns}
 \begin{cases}
    \begin{array}{l}
    \displaystyle u_t-\nu\Delta u+ (u \cdot \nabla ) u + \nabla p=0,\\
    \displaystyle\operatorname{div} u = 0,
    \end{array}&\quad\text{in}\quad\mathbb{R}^3\times(0,T]\\
    \ u(x,0)= u_0 ( x), &\quad \text{for any}\quad x \in \mathbb{R}^3.
  \end{cases}
\end{equation}
We naturally assume $u_0$ to be a solenoidal field.
The existence of smooth solution to (\ref{eq:ns}) is one of the famous millennium problems, posted by the Clay Mathematical Institute (C. L. Fefferman \cite{MR2238274}).

Leray \cite{MR1555394} established the existence of weak solutions to (\ref{eq:ns}). Hopf \cite{MR50423} continued on this and the weak solution today is known as the Leray-Hopf weak solution (see also \cite{MR0254401}):
\begin{theorem}[Leray-Hopf]
  Suppose $u_0\in L^2(\mathbb{R}^3)$ is solenoidal, then there exists a weak solution $u\in L^\infty(0,T;L^2(\mathbb{R}^3))\cap L^2(0,T;\dot{H}^1(\mathbb{R}^3))$ of (\ref{eq:ns}).
\end{theorem}

To solve the millennium problems, it suffices to establish a suitable integrability for the Leray-Hopf solution. The following spaces $L^{s,r}$ are often used to specify the integrability of $u$. Here, we are working on $\mathbb{R}^3$ or the case $n=3$.
\begin{equation}
L^{s,r}
    =\{u:\mathbb{R}^n\times (0,T)\rightarrow \mathbb{R}^n \,\,\mid \,\,\|u\|_{L^{s,r}}:=[\int^T_0 (\int_{\mathbb{R}^n} |u(x,t)|^sdx)^{\frac{r}{s}} dt]^{\frac{1}{r}}<\infty \}.
\end{equation}
Indeed, one has the following smoothness and uniqueness criterion.

\begin{theorem}
\label{thm:LadyProSer}
  Suppose $u_0\in L^2(\mathbb{R}^3)$, and $u$ is a Leray-Hopf solution of (\ref{eq:ns}). Moreover, suppose $u$ satisfies the following condition: for some $T=T_0>0$,
  \begin{equation}
  \label{cond:LadyProSer}
    u\in L^{s,r} \qquad\text{with}\quad \frac{3}{s}+\frac{2}{r}=1,\quad3\leqslant s\leqslant +\infty.
  \end{equation}
  Then $u$ is smooth in $\mathbb{R}^3\times (0,T_0]$ and any Leray-Hopf weak solutions of (\ref{eq:ns}) must coincide with $u$ in $\mathbb{R}^3\times (0,T_0]$.
\end{theorem}
Condition (\ref{cond:LadyProSer}) is now called the Serrin-Prodi-Ladyzhenskaya condition.
 For $3<s\leqslant +\infty$, the uniqueness was established by
 Serrin \cite{MR0150444} and Prodi \cite{MR126088} and the smoothness
 was proved by Ladyzhenskaya \cite{MR0236541}. L. Escauriaza, G. Seregin and
 V. \v{S}ver\'{a}k \cite{MR1992563} considered the limiting case $s=3$,
 they proved that any $L^\infty(0,T_0;L^3(\mathbb{R}^3))$ solution must be
 regular and can be extended further (see also the quantitative version by Terence Tao\cite{MR4337421}).

On the other hands, one can check that the Leray-Hopf weak solution must satisfies:
  \begin{equation}
  \label{cond:LerayHopf}
    u\in L^{s,r}\qquad\text{for}\quad \frac{3}{s}+\frac{2}{r}=\frac{3}{2}\qquad\text{and}\quad 2\leqslant r\leqslant +\infty.
  \end{equation}
There is a integrability gap between (\ref{cond:LadyProSer}) and (\ref{cond:LerayHopf}). As a step forward, Scheffer \cite{MR454426,MR510154,MR573611,MR676002} and Caffarelli-Khon-Nirenberg \cite{MR673830},established the partial regularity of suitable weak solutions to (\ref{eq:ns}). Lin \cite{MR1488514}  simplified the proof (See also \cite{MR1738171}). Indeed, the following result is obtained.
\begin{proposition}
\label{prop:Partial}
  There exists universal constants $\epsilon_0, C_0>0$, with the following property. Suppose $(u,p)$ is a suitable weak solution of (\ref{eq:ns})
\begin{equation}
  u_t+(u\cdot\nabla) u-\Delta u+\nabla p=0\qquad\text{in}\quad K_R(z_0,t_0):=(t_0-R^2,t_0)\times B_{R_0}(x_0).
\end{equation}
Let
\begin{equation}
  C(z_0,t_0,R):=R^{-2}\int_{K_R(z_0,t_0)}|u|^3\intd x\intd t\quad\text{and}\quad D(z_0,t_0R):=R^{-2}\int_{K_R(z_0,t_0)}|p-\bar{p}|^{\frac{3}{2}}\intd x\intd t.
\end{equation}
Assume $C(z_0,t_0,R)$ and $D(z_0,t_0,R)$ satisfy
\begin{equation}
  C(z_0,t_0,R)+D(z_0,t_0,R)\leqslant \epsilon_0.
\end{equation}
Then
\begin{equation}
  R|u|\leqslant C_0\qquad\text{in}\quad K_{R/2}(z_0,t_0),
\end{equation}
and $|\nabla^k u|$ is H\"{o}lder continuous in $K_{R/2}(z_0,t_0)$ for any $k\in\mathbb{N}$.
\end{proposition}

In 1972, Fabes, Jones and Riviere \cite{MR316915} established the short time existence of mild solution to (\ref{eq:ns}) satisfying condition (\ref{cond:LadyProSer}) for initial data $u_0\in L^s(\mathbb{R}^3)$ with $s>3$. Where the existing time $T_0$ depends on the initial data $u_0$. Moreover they obtained that if $u_0\in L^2(\mathbb{R}^3)\cap L^s(\mathbb{R}^3)$, then the mild solution is a Leray-Hopf weak solution, hence by Theorem~\ref{thm:LadyProSer} a smooth solution.

In 1989, Kato \cite{MR760047} considered the limiting case. For $u_0\in L^3(\mathbb{R}^3)$, the short time existence of mild solution. In his work, the existing time $T_0$ depends only on the $L^3$-norm of $u_0$, and that for sufficiently small $\|u_0\|_{L^3(\mathbb{R}^3)}$, it holds that $T_0=+\infty$. The following is a simple case that the small initial data implies the global existence:
\begin{proposition}
\label{prop:1}
  Suppose  $u_0$ satisfies:
  \begin{equation}
    \|u_0\|_{L^3(\mathbb{R}^3)}\leqslant C_\ast \nu,
  \end{equation}

  then
  problem (\ref{eq:ns}) has a  unique smooth solution,
\end{proposition}

%
%
%
%
%
%
%
By the well known energy inequality, one have the following:
\begin{proposition}
\label{prop:2}
For sufficiently large $T$, for example $T=C\|u_0\|_{L^2(\mathbb{R}^3)}^4\nu^{-4}$, there exists $t_0<T$ such that
\begin{equation}
   \|u(t_0)\|_{L^2(\mathbb{R}^3)}\|\nabla u(t_0)\|_{L^2(\mathbb{R}^3)}\leqslant C\nu^2.
\end{equation}
As a corollary of Proposition \ref{prop:1}, solution of (\ref{eq:ns}) is smooth for $t\geqslant t_0$.
\end{proposition}

In order to solve the millennium problem, we need to show global existence of
smooth solution for large initial data.
In \cite{MR2388660}, T. Hou and C. Li established the global regularity of solutions
to the incompressible axis-symmetric Navier-Stokes equations
in 3-d for some large initial data of special types.
Inspired by the existence of global smooth solutions for globally small initial data,
we aim at utilizing the local smallness and develop the global existence for
general compact supported smooth initial data.
An interesting counter part about the millennium problem can be found in
 \cite{hou2022potentially} by Thomas Hou and references there in.
They are seeking possible finite time blow-up of Leray-Hopf solutions
with certain regular initial data.
In this article, we derive the regularity of solutions to (\ref{eq:ns})
under local smallness assumption (\ref{eq:asu}).

Indeed, we define a localized $L^s$-norm of $u$ as:
\begin{equation}
  \|u\|_{L^s_R}=\sup \{\|u\|_{L^s({B_{R}(x)})} \mid x\in\mathbb{R}^n\}
\end{equation}

Our result states as follow:
\begin{theorem}
\label{thm:main}
  Suppose $u$ is a solution of (\ref{eq:ns}) satisfying
  \begin{equation}
  \label{eq:asu}
    u_0\in H^1(\mathbb{R}^3)\quad \int_0^T\|u(t)\|_{L^s_{R(t)}}^rdt<\infty\quad\text{for}\quad \frac{3}{s}+\frac{2}{r}=1\quad r<+\infty,
  \end{equation}
  where $R(t)$ is a positive function satisfying:
  \begin{equation}
  \label{eq:R}
    \int_0^T(R(t))^{-2}dt<\infty.
  \end{equation}
  Then $u\in L^\infty[0,T;C^\infty(\mathbb{R}^3)]$ with
  \begin{equation}
  \label{eq3}
\|\nabla u(t)\|_{L^2(\mathbb{R}^3)}\leqslant \|\nabla u(0)\|_{L^2(\mathbb{R}^3)}\exp\left\{\frac{C_1}{\nu^{r-1}}\int_0^t \|u(\tau)\|_{L^s_{R(\tau)}}^{r}d\tau +C_2\nu\int_{0}^tR(\tau)^{-2}d\tau\right\}.
  \end{equation}

\end{theorem}
\begin{remark}
Combining our result with Proposition \ref{prop:2}, we know that for given $u_0$, it suffices to require
(\ref{eq:asu}) and (\ref{eq:R}) for $T=C_*^{-2}\|u_0\|_{L^2(\mathbb{R}^3)}^4\nu^{-4}$ to derive the global smoothness: $u\in L^\infty[0,T;C^\infty(\mathbb{R}^3)]$.
\end{remark}
\begin{remark}
Our estimate appears to be similar to Theorem~\ref{thm:LadyProSer},
but only require a local integrability of $u$.
\end{remark}
%
The following is a special case of the Gargliado-Nirenberg inequalities on bounded domains (see \cite{MR4513000} by C. Li and K. Zhang).
\begin{lemma}
\label{lem:GN}
If $w \in H^2(Q)$ then
\begin{equation}
\int_{Q}|\nabla w-\overline{\nabla w}|^3dx\leqslant C\|w\|_{L^3(Q)}\|\nabla^2 w\|^2_{L^2(Q)}.
\end{equation}
 Here $\overline{\nabla w}$ is the average of $\nabla w$ over $Q$.
\end{lemma}
Next is the main estimate of this article.
It is interesting in its own and the method of derivation might be used to deal with other nonlinear problems.
\begin{lemma}[Main estimate]
\label{lem:est}
For $s$, there is an universal constant $C$, such that for any $\epsilon >0$:
\begin{equation}\label{eq15}
\left|\int_{\mathbb{R}^3}\frac{\partial u_k}{\partial x_i}\frac{\partial u_j}{\partial x_k}\frac{\partial u_j}{\partial x_i}dx\right|
\leqslant C_0\|u\|_{L^s_{\epsilon}} \left(\epsilon^{-\frac{3}{s}-1}\|\nabla u\|^2_{L^2(\mathbb{R}^3)}+\epsilon^{1-\frac{3}{s}}\|\nabla^2 u\|^2_{L^2(\mathbb{R}^3)}\right).
\end{equation}
\end{lemma}

\section{Proof of theorems}

\begin{proof}[Proof of Theorem \ref{thm:main}]
Without loss of generality, let $\nu=1$.
Define $$H(t)=\int_{\mathbb{R}^3} |\nabla u(x,t)|^2dx=\|\nabla u(t)\|^2_{L^2(\mathbb{R}^3)},$$ then we have the basic estimate:
\begin{equation}\label{eq5}
  \int_0^{t}H(\tau)d\tau=\int_0^{t}\int_{\mathbb{R}^3} |\nabla u(x,\tau)|^2dxd\tau =\frac{1}{2}(\|u_0\|^2_{L^2(\mathbb{R}^3)}-\|u(t)\|^2_{L^2(\mathbb{R}^3)}).
\end{equation}
Therefore
\begin{equation}
   \int_0^{T}H(t)dt=\int_{0}^T\|\nabla u(t)\|_{L^2(\mathbb{R}^3)}^2dt\leqslant\|u_0\|_{L^2(\mathbb{R}^3)}^2\qquad\text{and}\qquad \|u(t)\|_{L^2(\mathbb{R}^3)}\leqslant\|u_0\|_{L^2(\mathbb{R}^3)}.
\end{equation}

We can assume that $u(x,t) \in C^k((0,T], H^k(\mathbb{R}^3))$ with the help of classical regularity or approximation (approximate $u_0$
with $C_0^{\infty}(\mathbb{R}^3)$ functions) theory as long as we have a priori bounds on $\|u(\cdot, t)\|_{H^1(\mathbb{R}^3)}$, see \cite{MR1555394} and \cite{MR0150444}.

Next, we take the partial derivative $\frac{\partial}{\partial x_i}$ to equation (\ref{eq:ns}), take inner product with $\frac{\partial u}{\partial x_i}$, integral over $\mathbb{R}^3$ and then sum-up $i$ from 1 to 3 to get:
\begin{equation}\label{eq10}
  \frac{d}{dt}H(t)=\frac{d}{dt}\int_{\mathbb{R}^3} |\nabla u(x,t)|^2dx
  =-2\|\nabla^2 u(t)\|^2_{L^2(\mathbb{R}^3)}
  -2\sum_{i,j,k=1,2,3}\int_{\mathbb{R}^3}\left(\frac{\partial u_k}{\partial x_i}\frac{\partial u_j}{\partial x_k}\frac{\partial u_j}{\partial x_i}\right)_{(x,t)}dx.
\end{equation}

We apply the estimate in Lemma \ref{lem:est} to derive:
\begin{equation}\label{eq20}
  H'(t)\leqslant-2\|\nabla^2 u(t)\|^2_{L^2(\mathbb{R}^3)}+2C_0\|u(t)\|_{L^s_{\epsilon}} \left(\epsilon^{-\frac{3}{s}-1}\|\nabla u(t)\|^2_{L^2(\mathbb{R}^3)}+\epsilon^{1-\frac{3}{s}}\|\nabla^2 u(t)\|^2_{L^2(\mathbb{R}^3)}\right).
\end{equation}

Take
\begin{equation}
\epsilon=\min\left\{R(t),(C_0\|u(t)\|_{L^s_{R(t)}})^{-\frac{s}{s-3}}\right\}.
\end{equation}
Then
\begin{equation}
  \|u(t)\|_{L^s_{\epsilon}}\leqslant\|u(t)\|_{L^s_{R(t)}},\quad\text{since}\quad \epsilon\leqslant R(t).
\end{equation}
And therefore,
\begin{equation}
   2C_0\epsilon^{1-\frac{3}{s}}\|u(t)\|_{L^s_{\epsilon}}\leqslant 2C_0\epsilon^{1-\frac{3}{s}}\|u(t)\|_{L^s_{R(t)}}\leqslant 2.
\end{equation}

Then
\begin{equation}\label{eq30}
\begin{aligned}
H'(t)&\leqslant 2C_0\epsilon^{-1-\frac{3}{s}}\|u(t)\|_{L^s_{\epsilon}}H(t)\\
&\leqslant 2C_0\max\left\{(C_0\|u(t)\|_{L^s_{R(t)}})^{\frac{s+3}{s-3}}, R(t)^{-1-\frac{3}{s}}\right\}\|u(t)\|_{L^s_{R(t)}}H(t)\\
&\leqslant\left(2C_0^{\frac{2s}{s-3}}\|u(t)\|_{L^s_{R(t)}}^{\frac{2s}{s-3}} +2C_0\|u(t)\|_{L^s_{R(t)}}R(t)^{-1-\frac{3}{s}}\right)H(t)\\
&\leqslant\left(2C_0^{\frac{2s}{s-3}}\|u(t)\|_{L^s_{R(t)}}^{\frac{2s}{s-3}} +\frac{s-3}{2s}(2C_0)^{\frac{2s}{s-3}}\|u(t)\|_{L^s_{R(t)}}^{\frac{2s}{s-3}} +\frac{s+3}{2s}R(t)^{-2}\right)H(t)\\
&=\left(2C_1\|u(t)\|_{L^s_{R(t)}}^{\frac{2s}{s-3}}+2C_2R(t)^{-2}\right)H(t).
\end{aligned}
\end{equation}
Applying the Gronwall inequality, one obtains:
$$H(t)\leqslant H(0)\exp\left\{2C_1\int_0^t \|u(\tau)\|_{L^s_{R(\tau)}}^{\frac{2s}{s-3}}d\tau +2C_2\int_{0}^tR(\tau)^{-2}d\tau\right\}.$$
(\ref{eq3}) then follows.
\end{proof}

Now, we come back to prove the main estimate (\ref{eq15}).

\begin{proof}[Proof of Lemma~\ref{lem:est}]

We present our proof with the case $i=k=j$ and $u_i=u_j=u_k=w$, all other cases can be proved similarly. We demonstrate with the following:
\begin{equation}
\left|\int_{\mathbb{R}^3}(Dw)^3dx\right|
\leqslant C\max_{x\in\mathbb{R}^3}\|w\|_{L^s(Q_{\epsilon}(x))}\left(\epsilon^{-\frac{3}{s}-1}\|\nabla w\|^2_{L^2(\mathbb{R}^3)} +\epsilon^{1-\frac{3}{s}}\|\nabla^2 w\|^2_{L^2(\mathbb{R}^3)}\right).
\end{equation}

Here $D$ denotes a first order directional derivative. We decompose $\mathbb{R}^3$ into non-overlapping standard cubes $Q_{\epsilon/2}(x_i)$ of side length $\epsilon$.
By a little shifting of each edge of cubes $Q_{\epsilon/2}(x_i)$, one can make a re-decomposition as
\begin{equation}
  \mathbb{R}^n=\bigcup_{i}\tilde{Q}^i,
\end{equation}
such that for each $\tilde{Q}^i$ coming from $Q_{\epsilon/2}(x_i)$ the following inequality holds:
\begin{equation}
  \frac{1}{\epsilon^2}\int_{\partial \tilde{Q}^i}|w|d\sigma_x\leqslant \frac{C}{\epsilon^3}\int_{Q_\epsilon(x_i)}|w|dx.
\end{equation}

Now, let $$\overline{Dw}|_i=\frac{1}{|\tilde{Q}_i|}\int_{\tilde{Q}^i}Dwdx$$
be the average of $Dw$ on $\tilde{Q}^i$. We derive:
\begin{equation}\label{eq43}
\begin{aligned}
\int_{\mathbb{R}^3}(Dw)^3dx&=\sum_{i}\int_{\tilde{Q}^i}(Dw)^3dx\\
&=\sum_{i}\left[\int_{\tilde{Q}^i}(Dw-\overline{Dw}|_i)^3dx+3\overline{Dw}|_i\int_{\tilde{Q}^i}
(Dw-\overline{Dw}|_i)^2dx+|\tilde{Q}_i|\cdot\overline{Dw}|_i^3\right].
\end{aligned}
\end{equation}

We apply Lemma \ref{lem:GN} to estimate:
\begin{equation} \label{eq45}
  \begin{aligned}
  \int_{\tilde{Q}^i}\left|\nabla w-\overline{\nabla w}|_i\right|^3dx
  &\leqslant C\|w\|_{L^3(\tilde{Q}^i)}\|\nabla^2 w\|^2_{L^2(\tilde{Q}^i)} \\
  &\leqslant C \epsilon^{\frac{s-3}{s}}\|w\|_{L^s(\tilde{Q}^i)}\|\nabla^2 w\|^2_{L^2(\tilde{Q}^i)}.
\end{aligned}
\end{equation}

Therefore
\begin{equation}
\label{eq46}
\begin{aligned}
  \left|\sum_{i}\int_{\tilde{Q}^i}(Dw-\overline{Dw}|_i)^3dx\right|
  &\leqslant \sum_{i}\int_{\tilde{Q}^i}\left|\nabla w-\overline{\nabla w}|_i\right|^3dx \\
  &\leqslant C\epsilon^{\frac{s-3}{s}}\sup_{i}{\|w\|_{L^s(\tilde{Q}^i)}} \sum_{i} \|\nabla^2 w\|^2_{L^2(\tilde{Q}^i)}\\
  &=C\epsilon^{\frac{s-3}{s}}\sup_{i}{\|w\|_{L^s(\tilde{Q}^i)}} \|\nabla^2 w\|^2_{L^2(\mathbb{R}^3)}
  \end{aligned}
\end{equation}

On the other hand
\begin{equation}
  \begin{aligned}
    \left|3\overline{Dw}\int_{\tilde{Q}^i}(Dw-\overline{Dw}|_i)^2dx +|\tilde{Q}_i|\cdot\overline{Dw}|_i^3\right|
    &\leqslant C\epsilon^{-3}\|\nabla w\|^2_{L^2(\tilde{Q}^i)}\left|\int_{\tilde{Q}^i}Dwdx\right|\\
    &\leqslant C\epsilon^{-3} \|\nabla w\|^2_{L^2(\tilde{Q}^i)}\int_{\partial \tilde{Q}^i}|w|d\sigma_x\\
    &\leqslant C\epsilon^{-4} \|\nabla w\|^2_{L^2(\tilde{Q}^i)}\int_{Q_\epsilon(x_i)}|w|dx\\
    &\leqslant C\epsilon^{-\frac{3}{s}-1} \|\nabla w\|^2_{L^2(Q_\epsilon(x_i))}\|w\|_{L^s(Q_\epsilon(x_i))}.
  \end{aligned}
\end{equation}

Then
\begin{equation}
\label{eq47}
  \begin{aligned}
    \sum_i\left|3\overline{Dw}\int_{\tilde{Q}^i}(Dw-\overline{Dw})^2dx +|\tilde{Q}_i|\cdot\overline{Dw}^3\right|
    &\leqslant C\epsilon^{-\frac{3}{s}-1}\sup_i\|w\|_{L^s(Q_\epsilon(x_i))} \sum_i\|\nabla w\|^2_{L^2(Q_\epsilon(x_i))}\\
    &=C\epsilon^{-\frac{3}{s}-1}\sup_i\|w\|_{L^s(Q_\epsilon(x_i))}\|\nabla w\|^2_{L^2(\mathbb{R}^3)}.
  \end{aligned}
\end{equation}
Substituting (\ref{eq46}) and (\ref{eq47}) into (\ref{eq43}), note that
\begin{equation}
  \sup_i\|w\|_{L^s(\tilde{Q}^i)}\leqslant \sup_{x\in\mathbb{R}^3}\|w\|_{L^s(Q_{\epsilon}(x))}=\|w\|_{L^s_\epsilon},
\end{equation}
then we finish the proof.
\end{proof}




\section{Acknowledgement}
The authors are partially supported by National Natural Science Foundation of China (Grant Nos.12031012, 11831003).



\end{document}